\documentclass[11pt]{article}

\begin{document}

%%-I include macros files here-%%

%%dung le's macros

\newcommand{\newc}{\newcommand}

%%index

\renewcommand{\theequation}{\thesection.\arabic{equation}}
\newc{\eqnoset}{\setcounter{equation}{0}}
\newcommand{\myref}[2]{#1~\ref{#2}}

\newcommand{\mref}[1]{(\ref{#1})}
\newcommand{\reflemm}[1]{Lemma~\ref{#1}}
\newcommand{\refrem}[1]{Remark~\ref{#1}}
\newcommand{\reftheo}[1]{Theorem~\ref{#1}}
\newcommand{\refdef}[1]{Definition~\ref{#1}}
\newcommand{\refcoro}[1]{Corollary~\ref{#1}}
\newcommand{\refprop}[1]{Proposition~\ref{#1}}
\newcommand{\refsec}[1]{Section~\ref{#1}}
\newcommand{\refchap}[1]{Chapter~\ref{#1}}

%%environments
\newcommand{\beq}{\begin{equation}}
\newcommand{\eeq}{\end{equation}}
\newcommand{\beqno}[1]{\begin{equation}\label{#1}}

\newcommand{\barr}{\begin{array}}
\newcommand{\earr}{\end{array}}

\newc{\bearr}{\begin{eqnarray*}}
\newc{\eearr}{\end{eqnarray*}}

\newc{\bearrno}[1]{\begin{eqnarray}\label{#1}}
\newc{\eearrno}{\end{eqnarray}}

\newc{\non}{\nonumber}
\newc{\nol}{\nonumber\nl}

\newcommand{\bdes}{\begin{description}}
\newcommand{\edes}{\end{description}}
\newc{\benu}{\begin{enumerate}}
\newc{\eenu}{\end{enumerate}}
\newc{\btab}{\begin{tabular}}
\newc{\etab}{\end{tabular}}

%%\newtheorem{theorem}{Theorem}
%%\newtheorem{defi}[theorem]{Definition}
%%\newtheorem{lemma}{Lemma}[section]
%%\newtheorem{rem}[lemma]{Remark}
%%\newtheorem{exam}[lemma]{Example}
%%\newtheorem{propo}[theorem]{Proposition}
%%\newtheorem{corol}[theorem]{Corollary}

%%Unmark these for uniform indexing

\newtheorem{theorem}{Theorem}[section]
\newtheorem{defi}[theorem]{Definition}
\newtheorem{lemma}[theorem]{Lemma}
\newtheorem{rem}[theorem]{Remark}
\newtheorem{exam}[theorem]{Example}
\newtheorem{propo}[theorem]{Proposition}
\newtheorem{corol}[theorem]{Corollary}

\renewcommand{\thelemma}{\thesection.\arabic{lemma}}

\newcommand{\btheo}[1]{\begin{theorem}\label{#1}}
\newc{\brem}[1]{\begin{rem}\label{#1}\em}
\newc{\bexam}[1]{\begin{exam}\label{#1}\em}
\newc{\bdefi}[1]{\begin{defi}\label{#1}}
\newcommand{\blemm}[1]{\begin{lemma}\label{#1}}
\newcommand{\bprop}[1]{\begin{propo}\label{#1}}
\newcommand{\bcoro}[1]{\begin{corol}\label{#1}}
\newcommand{\etheo}{\end{theorem}}
\newcommand{\elemm}{\end{lemma}}
\newcommand{\eprop}{\end{propo}}
\newcommand{\ecoro}{\end{corol}}
\newc{\erem}{\end{rem}}
\newc{\eexam}{\end{exam}}
\newc{\edefi}{\end{defi}}

\newc{\rmk}[1]{{\bf REMARK #1: }}
\newc{\DN}[1]{{\bf DEFINITION #1: }}

\newcommand{\bproof}{{\bf Proof:~~}}
\newc{\eproof}{{\vrule height8pt width5pt depth0pt}\vspace{3mm}}

%%symbols

\newcommand{\rarr}{\rightarrow}
\newcommand{\Rarr}{\Rightarrow}
\newcommand{\tru}{\backslash}
\newc{\bfrac}[2]{\dspl{\frac{#1}{#2}}}

%%space

\newc{\nl}{\vspace{2mm}\\}
\newc{\nid}{\noindent}

%%operations

\newcommand{\oneon}[1]{\frac{1}{#1}}
\newcommand{\dspl}{\displaystyle}
\newc{\grad}{\nabla}
\newc{\Div}{\mbox{div}}
\newc{\pdt}[1]{\dspl{\frac{\partial{#1}}{\partial t}}}
\newc{\pdn}[1]{\dspl{\frac{\partial{#1}}{\partial \nu}}}
\newc{\pdNi}[1]{\dspl{\frac{\partial{#1}}{\partial \mathcal{N}_i}}}
\newc{\pD}[2]{\dspl{\frac{\partial{#1}}{\partial #2}}}
\newc{\dt}{\dspl{\frac{d}{dt}}}
\newc{\bdry}[1]{\mbox{$\partial #1$}}
\newc{\sgn}{\mbox{sign}}

\newc{\Hess}[1]{\frac{\partial^2 #1}{\pdh z_i \pdh z_j}}
\newc{\hess}[1]{\partial^2 #1/\pdh z_i \pdh z_j}

%%function spaces

\newcommand{\Coone}[1]{\mbox{$C^{1}_{0}(#1)$}}
\newcommand{\lspac}[2]{\mbox{$L^{#1}(#2)$}}
\newc{\hspac}[2]{\mbox{$C^{0,#1}(#2)$}}
\newc{\Hspac}[2]{\mbox{$C^{1,#1}(#2)$}}
\newc{\Hosp}{\mbox{$H^{1}_{0}$}}
\newcommand{\Lsp}[1]{\mbox{$L^{#1}(\Og)$}}
\newc{\hsp}{\Hosp(\Og)}

%%greeks

\newc{\ag}{\alpha}
\newc{\bg}{\beta}
\newc{\cg}{\gamma}\newc{\Cg}{\Gamma}
\newc{\dg}{\delta}\newc{\Dg}{\Delta}
\newc{\eg}{\varepsilon}
\newc{\zg}{\zeta}
\newc{\thg}{\theta}
\newc{\llg}{\lambda}\newc{\LLg}{\Lambda}
\newc{\kg}{\kappa}
\newc{\rg}{\rho}
\newc{\sg}{\sigma}\newc{\Sg}{\Sigma}
\newc{\tg}{\tau}
\newc{\fg}{\phi}\newc{\Fg}{\Phi}
\newc{\vfg}{\varphi}
\newc{\og}{\omega}\newc{\Og}{\Omega}
%\newc{\ng}{\eta}
\newc{\pdh}{\partial}

%%Integration and Sum

\newc{\ii}[1]{\int_{#1}}
\newc{\iidx}[2]{{\dspl\int_{#1}~#2~dx}}
\newc{\bii}[1]{{\dspl \ii{#1} }}
\newc{\biii}[2]{{\dspl \iii{#1}{#2} }}
\newc{\su}[2]{\sum_{#1}^{#2}}
\newc{\bsu}[2]{{\dspl \su{#1}{#2} }}

\newcommand{\iiomdx}[1]{{\dspl\int_{\Og}~ #1 ~dx}}
\newc{\biiom}[1]{{\dspl\int_{\bdrom}~ #1 ~d\sg}}
\newc{\io}[1]{{\dspl\int_{\Og}~ #1 ~dx}}
\newc{\bio}[1]{{\dspl\int_{\bdrom}~ #1 ~d\sg}}
\newc{\bsir}{\bsu{i=1}{r}}
\newc{\bsim}{\bsu{i=1}{m}}

\newc{\iibr}[2]{\iidx{\bprw{#1}}{#2}}
\newc{\Intbr}[1]{\iibr{R}{#1}}
\newc{\intbr}[1]{\iibr{\rg}{#1}}
\newc{\intt}[3]{\int_{#1}^{#2}\int_\Og~#3~dxdt}
%%\newc{\itQ}[2]{\dspl{\int\hspace{-2.5mm}\int_{#1}~#2~dxdt}}
%%\newc{\mitQ}[2]{\dspl{\rule[1mm]{4mm}{.3mm}\hspace{-5.3mm}\int\hspace{-2.5mm}\int_{#1}~#2~dxdt}}

\newc{\itQ}[2]{\dspl{\int\hspace{-2.5mm}\int_{#1}~#2~dz}}
\newc{\mitQ}[2]{\dspl{\rule[1mm]{4mm}{.3mm}\hspace{-5.3mm}\int\hspace{-2.5mm}\int_{#1}~#2~dz}}
\newc{\mitQQ}[3]{\dspl{\rule[1mm]{4mm}{.3mm}\hspace{-5.3mm}\int\hspace{-2.5mm}\int_{#1}~#2~#3}}

\newc{\mitx}[2]{\dspl{\rule[1mm]{3mm}{.3mm}\hspace{-4mm}\int_{#1}~#2~dx}}

\newc{\mitQq}[2]{\dspl{\rule[1mm]{4mm}{.3mm}\hspace{-5.3mm}\int\hspace{-2.5mm}\int_{#1}~#2~d\bar{z}}}
\newc{\itQq}[2]{\dspl{\int\hspace{-2.5mm}\int_{#1}~#2~d\bar{z}}}

\newc{\pder}[2]{\dspl{\frac{\partial #1}{\partial #2}}}

%%variables

\newc{\ui}{u_{i}}
\newcommand{\upl}{u^{+}}
\newcommand{\umn}{u^{-}}
\newcommand{\un}{\{ u_{n}\}}

\newcommand{\uo}{u_{0}}
\newc{\voi}{v_{i}^{0}}
\newc{\uoi}{u_{i}^{0}}
\newc{\vu}{\vec{u}}

\newc{\xo}{x_{0}}
\newc{\Br}{B_{R}}
\newc{\Bro}{\Br (\xo)}
\newc{\bdrom}{\bdry{\Og}}
\newc{\ogr}[1]{\Og_{#1}}
\newc{\Bxo}{B_{x_0}}

%%Additional macros
\newc{\inP}[2]{\|#1(\bullet,t)\|_#2\in\cP}
\newc{\cO}{{\mathcal O}}
\newc{\inO}[2]{\|#1(\bullet,t)\|_#2\in\cO}

\newc{\newl}{\\ &&}

%\newc{\bilhom}{\mbox{Bil}(\mbox{Hom}(\RR^n,\RR^N))}
\newc{\bilhom}{\mbox{Bil}(\mbox{Hom}(\RR^{nm},\RR^{nm}))}
\newc{\VV}[1]{{V(Q_{#1})}}

\newc{\ccA}{{\mathcal A}}
\newc{\ccB}{{\mathcal B}}
\newc{\ccC}{{\mathcal C}}
\newc{\ccD}{{\mathcal D}}
\newc{\ccE}{{\mathcal E}}
\newc{\ccH}{\mathcal{H}}
\newc{\ccF}{\mathcal{F}}
\newc{\ccI}{{\mathcal I}}
\newc{\ccJ}{{\mathcal J}}
\newc{\ccP}{{\mathcal P}}
\newc{\ccQ}{{\mathcal Q}}
\newc{\ccR}{{\mathcal R}}
\newc{\ccS}{{\mathcal S}}
\newc{\ccT}{{\mathcal T}}
\newc{\ccX}{{\mathcal X}}
\newc{\ccY}{{\mathcal Y}}
\newc{\ccZ}{{\mathcal Z}}

\newc{\bb}[1]{{\mathbf #1}}
\newc{\bbA}{{\mathbf A}}
%%Bracket
\newc{\myprod}[1]{\langle #1 \rangle}
\newc{\mypar}[1]{\left( #1 \right)}

%%Norms

\newc{\lspn}[2]{\mbox{$\| #1\|_{\Lsp{#2}}$}}
\newc{\Lpn}[2]{\mbox{$\| #1\|_{#2}$}}
\newc{\Hn}[1]{\mbox{$\| #1\|_{H^1(\Og)}$}}

%%Misc

\newc{\cyl}[1]{\og\times \{#1\}}
\newc{\cyll}{\og\times[0,1]}
\newc{\vx}[1]{v\cdot #1}
\newc{\vtx}[1]{v(t,x)\cdot #1}
\newc{\vn}{\vx{n}}

\newcommand{\RR}{{\rm I\kern -1.6pt{\rm R}}}

%\numberwithin{equation}{section}
%\newcommand{\eproof}{\vrule height5pt width3pt depth0pt}
%\newtheorem{exam}[lem]{Example}

\newenvironment{proof}{\noindent\textbf{Proof.}\ }
{\nopagebreak\hbox{ }\hfill$\Box$\bigskip}

%%Additional macros

\newc{\itQQ}[2]{\dspl{\int_{#1}#2\,dz}}
\newc{\mmitQQ}[2]{\dspl{\rule[1mm]{4mm}{.3mm}\hspace{-4.3mm}\int_{#1}~#2~dz}}
\newc{\MmitQQ}[2]{\dspl{\rule[1mm]{4mm}{.3mm}\hspace{-4.3mm}\int_{#1}~#2~d\mu}}

\newc{\MUmitQQ}[3]{\dspl{\rule[1mm]{4mm}{.3mm}\hspace{-4.3mm}\int_{#1}~#2~d#3}}
\newc{\MUitQQ}[3]{\dspl{\int_{#1}~#2~d#3}}

\vspace*{-.8in}
\begin{center} {\LARGE\em  Global Existence and Regularity Results for Strongly Coupled Nonregular Parabolic Systems via Iterative Methods.}

 \end{center}

\vspace{.1in}

\begin{center}

{\sc Dung Le}{\footnote {Department of Mathematics, University of
Texas at San
Antonio, One UTSA Circle, San Antonio, TX 78249. {\tt Email: Dung.Le@utsa.edu}
}}

\end{center}

\begin{abstract}
The global existence of classical solutions to strongly coupled parabolic systems is shown to be equivalent to the availability of an iterative scheme producing a sequence of solutions with uniform continuity in the BMO norms. Amann's results on global existence of classical solutions still hold under much weaker condition that their BMO norms do not blow up in finite time. The proof makes use of some new global and local weighted Gagliardo-Nirenberg inequalities involving BMO norms.
\end{abstract}

\vspace{.2in}

\section{Introduction}\label{introsec}\eqnoset

Among the long standing questions in the theory of strongly coupled parabolic systems and its applications are the global existence and regularity properties of their solutions. We consider in this paper the following system
\beqno{e1}\left\{\barr{ll} u_t=\Div(A(u)Du)+f(u)& (x,t)\in Q=\Og\times(0,T),\\u(x,0)=U_0(x)& x\in\Og\\\mbox{Boundary conditions for $u$ on $\partial \Og\times(0,T)$}. &\earr\right.\eeq

Here, $\Og$ is a bounded domain with smooth boundary $\partial \Og$ in $\RR^n$, $n>1$, and $u:\Og\to\RR^m$, $f:\RR^m\to\RR^m$ are vector valued functions. $A(u)$ is a full matrix $m\times m$. Thus, the above is a system of $m$ equations. The vector valued solution $u$ satisfies either Dirichlet or Neumann boundary condition on $\partial \Og\times(0,T)$.

The system \mref{e1} arises in many mathematical biology and ecology applications as well as in differential geometry theory. In the last few decades, papers concerning such strongly coupled parabolic systems usually assumed that the solutions under consideration were bounded, a very hard property to check as maximum principles had been unavailable for systems in general. In addition, past results usually relied on the following local existence result of Amann.

\btheo{Amthm} (\cite{Am1,Am2}) Suppose $\Og\subset\RR^n$, $n \ge2$, with $\partial\Og$
 being smooth. Assume that \mref{e1} is normally elliptic. Let $p_0 \in (n,\infty)$ and $U_0$ be  in $W^{1,p_0}(\Og
)$. Then there exists a maximal time $T_0\in(0,\infty]$ such that the system
\mref{e1} has a unique classical solution in $(0, T_0)$ with
$$u\in C([0,T_0),W^{1,p_0}(\Og)) \cap C^{1,2}((0, T_0)\times\bar{\Og})$$
Moreover, if $T_0<\infty$ then \beqno{blowup}\lim_{t\to T^-_0}\|u(\cdot,t)\|_{W^{1,p_0}(\Og)}=\infty.\eeq\etheo

The proof of the above result worked directly with the system and based on semigroup  and interpolation of functional spaces theories. We refer the readers to \cite{Am1} for the definition of normal ellipticity. The checking of \mref{blowup} is the most difficult one as known techniques for the regularity of solutions to scalar equations could not be extended to systems and counterexamples were available.  

In this paper, we propose a different approach using iterative techniques and depart from the boundedness assumptions. Namely, we consider the following schemes
$$(u_k)_t=\Div(A(u_{k-1})Du_k)+f(u_{k-1}, Du_k)\quad k\ge1.$$ 

Under very weak assumptions on the uniform boundedness and continuity of the BMO norms of the solutions to the above systems, we will show the global existence of a classical solution to \mref{e1}. Thus, global existence and regularity of solutions are established at once. Furthermore, without the boundedness assumptions the systems are no longer regular elliptic, we will only assume that the matrix $A(u)$ in \mref{e1} is {\em uniformly elliptic}.

We also improve \reftheo{Amthm} by replacing the condition \mref{blowup} with a weaker ones using the BMO or $W^{1,p_0}$ norms of $u$ with $p_0=n$. In a forthcoming work, we will show that the results in this paper can apply to a class of generalized Shigesada-Kawasaki-Teramoto models (\cite{SKT}) consisting of more than 2 equations. Namely, we will establish the global existence of classical solutions to the following system
\beqno{SKTgen} u_t = \Delta(P(u)) + f(u), \eeq where $P(u),f(u)$ are vector valued functions whose components have quadratic (or even polynomial) growth in $u$.

In the proof we make use of some new local Gagliardo-Nirenberg inequalities involving BMO norms and weights in $A_p$ classes. These are the generalizations of the inequalities by Strzelecki and Rivi\`ere in \cite{SR}. 

\section{Preliminaries and Main Results}\label{premsec}\eqnoset

Throughout this paper $\Og$ is a bounded domain with smooth boundary in $\RR^n$, $n>1$. To describe our assumptions we recall the definitions of BMO spaces and $A_p$ classes.

For any locally integrable vector valued function $u\in L^1_{loc}(\Og,\RR^m)$ and measurable set $B\subset\Og$ with its Lebesgue measure $|B|\ne0$, we denote $$\mitx{B}{u} = \frac{1}{|B|}\iidx{B}{u}.$$

For a smooth function $u$ defined on $\Og\times(0,T)$, $T>0$, its temporal and spatial derivatives are denoted by $u_t,Du$ respectively. If $A$ is a function in $u$ then we also abbreviate $\frac{\partial A}{\partial u}$ by $A_u$.

We will frequently work with a ball $B(z,R)$ centered at $z\in\RR^n$ with radius $R$. If $z$ is understood we will write $B_R$ for $B(z,R)$ 
For any $x_0\in\Og$ and $R>0$, we also denote $\Og(x_0,R)=\Og\cap B(x_0,R)$. 
For any locally measurable function $u$ on $\Og$ and $x_0\in\Og$ and $R>0$ we write $u_{x_0,R}$ for the average of $u$ over $\Og(x_0,R)$. If $x_0$ is understood, we simply write $u_R$ for $u_{x_0,R}$.

We say that a locally integrable vector valued function $u:\Og\to\RR^m$ is BMO (Bounded Mean Oscillation) if the seminorm $$[u]_{BMO(\Og)} = \sup_{x_0\in \Og, R>0}\mitx{\Og(x_0,R)}{\left|u-u_{\Og(x_0,R)}\right|} <\infty,$$ where the supremum is taken over all balls $B\subset \Og$. The space $BMO(\Og)$ is the Banach space of BMO functions on $\Og$ with norm $$\|u\|_{BMO(\Og)}=\|u\|_{L^1(\Og)} + [u]_{BMO(\Og)}.$$

We recall the following well known fact (e.g., see \cite{JN}).
\blemm{bmojn} 
If $\llg:\RR^m\to\RR$ is H\"older function then $\llg(u)$ is BMO if $u$ is BMO.\elemm

There is a connection between BMO functions and the so called $A_\cg$ weights, which  are defined as follows. Let  $\Psi$ be a measurable nonnegative function on $\Og$ and $\cg>1$. We say that $\Psi$ belongs to the class $A_\cg$ or $\Psi$ is an $A_\cg$ weight if the quantity
\beqno{aweight} [\Psi]_{\cg} = \sup_{B(x,R)\subset\Og} \left(\mitx{B(x,R)}{\Psi}\right) \left(\mitx{B(x,R)}{\Psi^{1-\cg'}}\right)^{\cg-1}< \infty.\eeq
Here, $\cg'=\cg/(\cg-1)$. The $A_\infty$ class is defined by $A_\infty=\cup_{\cg>1}A_\cg$. For more details on these classes we refer the readers to \cite{FPW,OP,st}.

We also recall the following result from \cite{JN} on the connection between BMO functions and weights.
\blemm{abmoweight} Let $\Psi$ be a positive function and $\mu$ is a nonnegative Radon measure on $\Og$ such that $\Psi,\Psi^{-1}$ are in $BMO(\Og)$ then $\Psi$ belongs to $\cap_{\cg>1}A_\cg$ and $[\Psi]_{\cg}$ is bounded by a constant depending on $[\Psi]_{BMO(\Og)}$ and $[\Psi^{-1}]_{BMO(\Og)}$. \elemm

We also recall the definition of the  Campanato spaces ${\cal L}^{p,\cg}(\Og,\RR^m)$. For any $p\ge1$ and $\cg>0$ and $u\in L^p(\Og,\RR^m)$, we define 
$$[u]_{p,\cg}^p = \sup_{x_0\in\Og, \rg>0} \rg^{-\cg}\iidx{\Og(x_0,\rg)}{|u-u_{x_0,\rg}|^p}.$$

Then ${\cal L}^{p,\cg}(\Og,\RR^m)$ is the Banach space of such functions with finite  norm $$\|u\|_{p,\cg} = \|u\|_{p} + [u]_{p,\cg}.$$ 

Clearly, ${\cal L}^{p,n}(\Og,\RR^m)=BMO(\Og,\RR^m)$. Moreover, it is well known that (\cite[Theorem 2.9, p.52]{Gius})  ${\cal L}^{p,\cg}(\Og,\RR^m)$ is isomorphic to $C^{0,\ag}(\Og)$ if $\ag=\frac{\cg-n}{p}>0$.

As usual, $W^{1,p}(\Og,\RR^m)$, $p\ge1$, will denote the standard Sobolev spaces whose elements are vector valued functions $u\,:\,\Og\to \RR^m$ with finite norm $$\|u\|_{W^{1,p}(\Og,\RR^m)} = \|u\|_{L^p(\Og)} + \|Du\|_{L^p(\Og)},$$ where $Du$ is the the derivative of $u$.

We now state our structural conditions on the system \mref{e1}.

\bdes \item[A.1)] (Uniform ellipticity) There are positive constants $C,\llg_0$ and smooth functions $\llg(u),\LLg(u)$ such that $\llg(u)\ge \llg_0$ and $\LLg(u) \le C\llg(u)$ for all $u\in\RR^m$ and $$\llg(u)|\xi|^2 \le \myprod{A(u)\xi,\xi} \le \LLg(u)|\xi|^2\quad \forall u\in\RR^m,\,\xi\in\RR^{nm}.$$

\item[A.2)] Assume that $A\in C^1(\RR^m)$. Let $\Fg_0,\Fg$ be defined as $$\Fg_0(u) = \llg^\frac12(u)\mbox{ and } \Fg(u)= \frac{|A_u(u)|}{\llg^\frac12(u)}\quad u\in\RR^m.$$ Assume that the quantities
\beqno{CFg}
k_1:=\sup_{u\in\RR^m}\frac{|\Fg_u|}{\Fg},\, k_2:=\sup_{u\in\RR^m}\frac{\Fg}{\Fg_0}\eeq are finite.
\edes

\bdes
\item[A.3)](Weights) If $u\in BMO(\Og)$ then $\Fg^\frac23(u)$ belongs to the $A_\frac43$ class and the quantity $[\Fg(u)^\frac23]_{\frac43}$ can be controlled by the norm $\|u\|_{BMO(\Og)}$.
\edes

In applications, if $A(u)$ has a polynomial growth in $u$ then we can assume that $|A(u)|\sim \llg(u)$ and $|A_u|\sim|\llg_u|$. In this case
$$\Fg\sim \frac{|\llg_u|}{\llg^\frac12},\,|\Fg_u|\sim \frac{|\llg_{uu}|}{\llg^\frac12}+\frac{|\llg_u|^2}{\llg^\frac32}.$$ Therefore, if $\sup_u\frac{|\llg_u|}{\llg}$ and $\sup_u\frac{|\llg_{uu}|}{|\llg_u|}$ are bounded then \mref{CFg} of A.2) is verified. It is clear that this is the case if $\llg(u)$ is a polynomial in $|u|$, say $\llg(u)\sim (1+|u|)^k$ for some $k\ge0$. Similarly, concerning A.3), we see that $\Fg^\ag(u)\sim(1+|u|)^{\ag(k/2-1)}$. Thus, if $2\le k<5$ then $\ag(k/2-1)\le1$ for some $\ag>2/3$. In this case, $\Fg^\ag(u)$ is H\"older in $u$ so that it is BMO if $u$ is. Of course, $\Fg^{-1}(u)$ is bounded so that it is also BMO. \reflemm{abmoweight} shows that $\Fg^\ag(u)$ is then an $A_\cg$ weight for all $\cg>1$ and A.3) is satisfied. 

On the other hand, it is clear from the definition of weights that $[\Fg^\ag]_{\bg+1} = [\Fg^{-\frac{\ag}{\bg}}]_\frac{\bg+1}{\bg}^\bg$ so that if $1<k<2$ then $\Fg^{\frac{\ag}{\bg}}$ is bounded from above and $\Fg^{-\frac{\ag}{\bg}}\sim (1+|u|)^{\frac\ag\bg(1-k/2)}$. Thus, we can find $\ag>2/3$ and $\bg<1/3$ such that $\Fg^{-\frac{\ag}{\bg}}$ is H\"older in $u$ and therefore BMO if $u$ is. Again, this gives that $\Fg^{-\frac{\ag}{\bg}}$ is a weight and belongs to $\cap_{\cg>1}A_\cg$.

Concerning global existence of classical solutions, we also assume that the ellipticity constants $\llg,\LLg$ in A.1) are not too far apart.
\bdes
\item[R)] (The ratio condition) There is $\dg\in[0,1)$ such that \beqno{ratio} \frac{n-2}{n} = \dg \sup_{u\in\RR^m}\frac{\llg(u)}{\LLg(u)}.\eeq
\edes

One should note that there are examples in \cite{sd} of blow up solutions to \mref{e1} if the condition R) is violated.

We assume the following growth conditions on the nonlinearity $f$.

\bdes \item[F)] There are positive constants $C,b$ such that for any vector valued functions $u\in C^1(\Og,\RR^m)$ and $p\in C^1(\Og,\RR^{mn})$ \beqno{fcond2a} |f(u,p)| \le C|p| + C|u|^b+C.\eeq
\beqno{fcond3a} |Df(u,p)| \le C|Dp| + C|u|^{b-1}|Du|+C.\eeq
\edes

To solve \mref{e1}, we can make use of the following iterative scheme. We start with any smooth vector valued function $u_0$ on $Q$ and define a sequence $\{u_k\}$ of solutions to the following linear systems
\beqno{ga1} (u_k)_t=\Div(A(u_{k-1})Du_k)+f(u_{k-1}, Du_k)\quad k\ge1.\eeq 

The initial and boundary conditions for the above systems are those of $u$ in \mref{e1}. Note that the global existence of the strong solutions to the above systems is not generally available by standard theories (see \cite{AF}) because of the presence of $Du_k$ in $f$. However, this is the case if we assume R) and the linear growth of $f$ in $Du$ of F) and make use of the results in \cite{sd}.

Concerning the approximation sequence $\{u_k\}$, we assume the following uniform bound and continuity of their BMO norms.

\bdes\item[V)] Let $\{u_k\}$ and $\Fg$ be defined by \mref{ga1} and A.2). There exists a continuous function $C$ on $(0,\infty)$ such that for any $T>0$ $$\|u_0(\cdot,t)\|_{C^1(\Og)},\, \|u_k(\cdot,t)\|_{BMO(\Og)} \le C(T) \quad \forall t\in(0,T),\, k=1,2,\ldots$$ Moreover, for any $\eg>0$ and $(x,t)\in Q$, there exists $R=R(\eg,T)>0$ such that \beqno{vmouni} \|u_k(\cdot,t)\|_{BMO(B_R(x))}<\eg\quad \forall t\in(0,T),\,k=0,1,\ldots\eeq

In addition, for all $(x,t)\in Q$ and integer $k\ge1$ we assume that
\beqno{UF}\mbox{either }\frac{\Fg(u_{k}(x,t))}{\Fg(u_{k-1}(x,t))}\le C(T) \mbox{ or } 
\frac{\llg(u_{k}(x,t))}{\llg(u_{k-1}(x,t))}\le C(T).\eeq 
\edes

The uniform boundedness assumption on the BMO norm of $u_k$ is of course much weaker than the $L^\infty$ boundedness assumptions in literature. Moreover, the uniform continuity assumption \mref{vmouni} on the BMO norms is somehow necessary for the regularity of the limit solution $u$.

The assumption \mref{UF} seems to be technical at first glance but it is clearly necessary if we would like to produce a the sequence $\{u_k\}$ that converges in $L^\infty(Q)$ to a solution of \mref{e1}. 
In paticular, if $\llg(u)$ behaves like $\llg(u)\sim (1+|u|)^k$ for some $k\ge0$ then, as discussed earlier, we see that $\Fg(u)\sim(1+|u|)^{(k/2-1)}$. Thus, $$\frac{\Fg(u_{k}(x,t))}{\Fg(u_{k-1}(x,t))} \sim \left(\frac{1+|u_{k}(x,t)|}{1+|u_{k-1}(x,t)|}\right)^{k/2-1},\,  
\frac{\llg(u_{k}(x,t))}{\llg(u_{k-1}(x,t))}\sim \left(\frac{1+|u_{k}(x,t)|}{1+|u_{k-1}(x,t)|}\right)^{k}.$$
Hence, if $k\in[0,2]$ then \mref{UF} is clearly verified. The generalized Shigesada-Kawasaki-Teramoto model \mref{SKTgen} clearly is a typical example.

We then have our main result of this paper as follows.
\btheo{GNthm1} Assume A.1)-A.3), R), F) and V). Then the sequence $\{u_k\}$ has a subsequence that converges strongly to a classical solution $u$ in $C^0(\Og\times(0,\infty))$.
\etheo

A simple consequence of the above theorem is the following.
\bcoro{GNthm1z} Assume A.1)-A.3), R), F) and V). Then the system \mref{e1} has a classical solution $u$ that exists globally on $\Og\times(0,\infty)$ if and only if there is $u_0$ and a sequence $\{u_k\}$  satisfying V) such that $\{u_k\}$ has a subsequence that converges weakly to $u$ in $L^2(\Og\times(0,T))$  for each $T>0$.
\ecoro

The necessary part is trivial as we can take $u_k=u$ for all $k\ge 0$, the solution sequence is then a constant one. The sufficient part comes from \reftheo{GNthm1} and the uniqueness of weak limits in $L^2(\Og\times(0,T))$.

In the next theorem we discuss the global existence of classical solutions when their local existence can be achieved by other methods (e.g., \reftheo{Amthm}). 

\btheo{Amthm1}  Assume A.1)-A.3), R) and F).  Let $p_0 \in (n,\infty)$ and $U_0$ be  in $W^{1,p_0}(\Og
)$. Suppose that $T_0\in(0,\infty]$ is the maximal existence time for a classical solution $$u\in C([0,T_0),W^{1,p_0}(\Og)) \cap C^{1,2}((0, T_0)\times\bar{\Og})$$ for the system
\mref{e1}. 
Suppose that there is a function $C$ in $C^0((0,T_0])$ such that $$ \|u(\cdot,t)\|_{BMO(\Og)} \le C(t) \quad \forall t\in(0,T_0).$$
Moreover, for any $\eg>0$ and $(x,t)\in Q$, there exists $R=R(\eg)>0$ such that \beqno{vmouni1} \|u(\cdot,t)\|_{BMO(B_R(x))}<\eg\quad \forall t\in(0,T_0).\eeq

Then $T_0=\infty$.
\etheo

By Poincar\'e's inequality, $$\mitx{B_R}{|u-u_R|^n}\le C(n)\iidx{B_R}{|Du|^n},$$
it is easy to see that if $Du(\cdot,t)\in L^n(\Og)$ then $u(\cdot,t)$ is BMO and $\|u(\cdot,t)\|_{BMO(B_R)}$ is small if $R$ is small. Therefore, as a simple consequence of the above theorem, we have the following improvisation of \reftheo{Amthm}.

\bcoro{Amz} In addition to the assumptions of \reftheo{Amthm}, we assume R). Then there exists a maximal time $T_0\in(0,\infty]$ such that the system
\mref{e1} has a unique classical solution in $(0, T_0)$ with
$$u\in C([0,T_0),W^{1,p_0}(\Og)) \cap C^{1,2}((0, T_0)\times\bar{\Og})$$
Moreover, if $T_0<\infty$ then \beqno{blowup1}\lim_{t\to T^-_0}\|u(\cdot,t)\|_{W^{1,n}(\Og)}=\infty.\eeq\ecoro

Again, we remark that if the condition R) is violated then there are counterexamples for finite time blow up solutions to \mref{e1}.

\section{Weighted Gagliardo-Nirenberg inequalities}\label{GNineq}\eqnoset
In this section we will establish global and local weighted Gagliardo-Nirenberg  interpolation inequalities which allow us to control the $L^{2p+2}$ norm of $Du_k$ in the proof of our main theorems.

Since the inequalities can be useful for other applications, we will prove them under a set of the following independent assumptions.

\bdes\item[GN.1)]Let $\Fg_0,\Fg$ be positive functions on $\RR^m$ with $\Fg\in C^1(\RR^m)$. Assume that the following quantities are finite
\beqno{CFg1}
k_1:=\sup_{u\in\RR^m}\frac{|\Fg_u|}{\Fg}, \, k_2:=\sup_{u\in\RR^m}\frac{\Fg}{\Fg_0}.\eeq

\item[GN.2)] For some $p\ge1$ suppose that $\Fg(u)^\frac{2}{p+2}$ belongs to the $A_{\frac{p}{p+2}+1}$ class if $u$ is BMO and the quantity $[\Fg^\frac{2}{p+2}]_{\frac{p}{p+2}+1}$ can be controlled by the norm $\|u\|_{BMO(\Og)}$.
\edes

In  the rest of this paper we will slightly abuse our notations and write the dot product $\myprod{u,v}$ as $uv$ for any two vectors $u,v$ because its meaning should be clear in the context. Similarly, when there is no ambiguity $C$ will denote a universal constant that can change from line to line in our argument. Furthermore, $C(\cdots)$ is used to denote quantities which are bounded in terms of theirs parameters. 

\blemm{GNlobal} Assume GN.1) and GN.2). Let $u,U:\Og\to \RR^m$ be vector valued functions with $U\in C^2(\Og)$, $u\in C^1(\Og)$. Suppose further that either $U$ or $\Fg^2(u)\frac{\partial U}{\partial \nu}$ vanish on the boundary $\partial \Og$  of $\Og$. We set \beqno{Idef} I_1:=\iidx{\Og}{\Fg^2(u)|DU|^{2p+2}},\,\hat{I}_1:=\iidx{\Og}{\Fg^2(u)|Du|^{2p+2}},\eeq and \beqno{Idef2}
I_2:=\iidx{\Og}{\Fg_0^2(u)|DU|^{2p-2}|D^2U|^2}.\eeq 

Then there is a constant $C_\Fg$ depending on $[\Fg^\frac{2}{p+2}(u)]_{\frac{p}{p+2}+1}$ for which \beqno{GNglobalest}I_1 \le  C_\Fg\|U\|_{BMO(\Og)}
\left[k_1(I_1+\hat{I}_1) + k_2I_1^\frac12 I_2^\frac12\right].\eeq
\elemm

Before we go to the proof of this lemma, let us recall the following facts from Harmonic Analysis.
We first recall the definition of the maximal function of a function $F\in L^1_{loc}(\Og)$
$$ M(F)(x) = \sup_\eg\{\iidx{B_\eg(x)}{F(x)}\,:\, \eg>0 \mbox{ and } B_\eg(x)\subset\Og\}.$$

We also note here the Hardy-Littlewood theorem for any $F\in L^q(\Og)$
\beqno{HL}\iidx{\Og}{M(F)^q} \le C(q)\iidx{\Og}{F^q}, \quad q>1.\eeq

In addition,  let us recall the definition of the Hardy space ${\cal H}^{1}$. For any $y\in\Og$ and $\eg>0$. Let $\fg$ be any function in $C^1_0(B_1(y))$ with $|D\fg|\le C_1$. Let $\fg_\eg(x)=\eg^{-n}\fg(\frac{x}{\eg})$ (then $|D\fg_\eg|\le C_1\eg^{-1}$). From \cite{S}, a function $g$ is in ${\cal H}^{1}(\Og)$ if $$\sup_{\eg>0}g*\fg_\eg\in L^1(\Og) \mbox{ and } \|g\|_{{\cal H}^1} = \|g\|_{L^1(\Og)}+\|\sup_{\eg>0}g*\fg_\eg\|_{L^1(\Og)}.$$

Furthermore, concerning the $A_\cg$ classes, it is well known (e.g. \cite[Corollary 2.5]{OP}) that if $w$ belongs to the $A_\cg$ class for some $\cg>1$ then there are positive $\eg,\nu$ depending on $[w]_\cg$ such that $w^{1+\eg}$ belongs to the $A_{\cg-\nu}$ class and $[w^{1+\eg}]_{\cg-\nu}$ can be controlled by $[w]_\cg$. Hence,
by GN.2) there are positive constants $\ag,\bg$ depending on the quantity $[\Fg^\frac{2}{p+2}]_{\frac{p}{p+2}+1}$ such that \beqno{agbg}\ag>2/(2+p), \,\bg<p/(p+2) \mbox{ and } [\Fg^{\ag}]_{\bg+1}\le C([\Fg^\frac{2}{p+2}]_{\frac{p}{p+2}+1}).\eeq

A simple use of H\"older's inequality also gives \beqno{wagbg} [w^\dg]_\cg\le [w]_\cg^\dg \quad \forall \dg\in(0,1).\eeq

\bproof Integrating by parts, we have
\beqno{ibp}\iidx{\Og}{\Fg^2(u)|DU|^{2p+2}} = -\iidx{\Og}{U\Div(\Fg^2(u)|DU|^{2p}DU)}.\eeq

We will show that  $g=\Div(\Fg^2(u)|DU|^{2p}DU)$ belongs to the Hardy space ${\cal H}^{1}$ and 
\beqno{gh1est}\|g\|_{{\cal H}^1} = \iidx{\Og}{\sup_\eg|g*\fg_\eg|} \le C([\Fg^\frac{2}{p+2}(u)]_{\frac{p}{p+2}+1})
\left[k_1(I_1+\hat{I}_1) + k_2I_1^\frac12 I_2^\frac12\right].\eeq 
Once this is established, \mref{ibp} and the duality of the BMO and Hardy spaces give \mref{GNglobalest}.

We write $g=g_1+g_2$ with $g_i=\Div V_i$, setting
$$V_1= \Fg(u) |DU|^{p+1}\left(\Fg(u) |DU|^{p-1}DU-\mitx{B_\eg}{\Fg(u) |DU|^{p-1}DU}\right),$$ and
$$V_2= \Fg(u) |DU|^{p+1}\mitx{B_\eg}{\Fg(u) |DU|^{p-1}DU}.$$

Let us consider $g_1$ first and define $h=\Fg(u) |DU|^{p-1}DU$. For any $y\in\Og$ and $B_\eg=B_\eg(y)\subset \Og$, we use integration by parts, the property of $\fg_\eg$ and then H\"older's inequality for any $s>1$ to have the following
$$\barr{lll}|g_1*\fg_\eg|&=&\left|\iidx{B_\eg}{D\fg(\frac{x-y}{\eg})(h-h_{B_\eg})\Fg(u)|DU|^{p+1}}\right|\\&\le&\frac{C_1}{\eg}
\left|\mitx{B_\eg}{|h-h_{B_\eg}|\Fg(u)|DU|^{p+1}}\right|
\\&\le&\frac{C_1}{\eg}
\left(\mitx{B_\eg}{|h-h_{B_\eg}|^s}\right)^\frac1s
\left(\mitx{B_\eg}{\Fg^{s'}(u)|DU|^{(p+1)s'}}\right)^\frac1{s'}.
\earr$$ 

By Sobolev-Poincar\'e's inequality, with $s_*=n s/(n +s)$, we have the following estimate in noting that $|Dh| \le |D\Fg||DU|^p + p\Fg|DU|^{p-1}|D^2U|$
\beqno{dhest}\barr{ll}\lefteqn{\frac{C_1}{\eg}
\left(\mitx{B_\eg}{|h-h_{B_\eg}|^s}\right)^\frac1s \le C
\left(\mitx{B_\eg}{|Dh|^{s_*}}\right)^\frac1{s_*}}\hspace{1cm}&\\& \le
C\left[\mitx{B_\eg}{|D\Fg|^{s_*}|DU|^{ps_*}}+
\mitx{B_\eg}{\Fg^{s_*}|DU|^{(p-1)s_*}|D^2U|^{s_*}}\right]^\frac1{s_*}
\earr\eeq Take $s=2n/(n-1)$ then $s_*=s'=2n/(n+1)<2$. From the definitions of $k_1=\sup_{u\in\RR^m}\frac{|\Fg_u|}{\Fg}$ and $k_2=\sup_{u\in\RR^m}\frac{\Fg}{\Fg_0}$ in GN.1) we have
$$\left(\mitx{B_\eg}{|D\Fg|^{s_*}|DU|^{ps_*}}\right)^\frac1{s_*} \le k_1
\left(\mitx{B_\eg}{\Fg|DU|^{ps_*}|Du|^{s_*}}\right)^\frac1{s_*},$$ and
$$\left(\mitx{B_\eg}{\Fg^{s_*}|DU|^{(p-1)s_*}|D^2U|^{s_*}}\right)^\frac1{s_*}\le
k_2 \left(\mitx{B_\eg}{\Fg_0^{s_*}|DU|^{(p-1)s_*}|D^2U|^{s_*}}\right)^\frac1{s_*}.$$

Using the above estimates in \mref{dhest}, we get
\beqno{hest}\frac{C_1}{\eg}
\left(\mitx{B_\eg}{|h-h_{B_\eg}|^s}\right)^\frac1s \le C\left[k_1\Psi_1+k_2\Psi_2\right],\eeq where $$\Psi_1(y)=\left(M(\Fg^{s_*}|DU|^{ps_*}|Du|^{s_*})(y)\right)^\frac1{s_*}, \quad  \Psi_2(y)=\left(M(\Fg_0^{s_*}|DU|^{(p-1)s_*}|D^2U|^{s_*})(y)\right)^\frac1{s_*}.$$

Putting these estimates together we thus have \beqno{g1a}\sup_{\eg>0,}|g_1*\fg_\eg| \le C
\Psi_1\left[k_1\Psi_1+k_2\Psi_2\right].\eeq 

 Because $2>2n/(n+1)=s_*$, we can use Young's inequality and then \mref{HL} to get
$$\barr{lll}\left(\iidx{\Og}{\Psi_1^2}\right)^\frac12 &\le& \|M(|\Fg|^{s_*}|DU|^{(p+1)s_*})\|_{\frac{2}{s_*}}^\frac{1}{s_*}+\|M(|\Fg|^{s_*}|Du|^{(p+1)s_*})\|_{\frac{2}{s_*}}^\frac{1}{s_*}\\&\le& C(\|\Fg^2|DU|^{2(p+1)}\|_2+\|\Fg^2|Du|^{2(p+1)}\|_2),\earr$$ and
$$\left(\iidx{\Og}{\Psi_2^{2}}\right)^\frac1{2} = \left\|M\left(\Fg_0^{s_*}|DU|^{(p-1)s_*}|D^2u|^{s_*}\right)\right\|_{\frac{2}{s_*}}^\frac{1}{s_*}\le C\|\Fg_0^2|DU|^{p-1}|D^2U|\|_{2}.$$

Therefore, by Holder's inequality and the above estimates \beqno{g1est}\iidx{\Og}{\sup_\eg|g_1*\fg_\eg|} \le C
\left[k_1(I_1+\hat{I}_1) + k_2I_1^\frac12 I_2^\frac12\right].\eeq

We now turn to $g_2$ and note that  $|\Div V_2| \le C(J_1+J_2)$ for some constant $C$ and $$J_1:= \sup\frac{|\Fg_u|}{\Fg}\Fg|DU|^{p+1}|Du|J_3, \quad  J_2:= \Fg|DU|^p|D^2U|J_3 \mbox{ with } J_3:=\left|\mitx{B_\eg}{\Fg |DU|^{p}}\right|.$$

We consider $J_1$. For any $r>1/(p+1)$ we denote $r^*=1-\frac{1}{r(p+1)}$. We also write $F=\Fg|DU|^{p+1}$ and $\hat{F}=\Fg|Du|^{p+1}$. For  any $s>0$, H\"older's inequality  yields  $$\left(\mitx{B_\eg}{\Fg^{\frac{-s}{(p+1)}}\Fg^{\frac{s}{p+1}}|Du|^s}\right)^\frac1s \le \left(\mitx{B_\eg}{\Fg^{\frac{-s}{r^*(p+1)}}}\right)^{\frac{r^*}{s}}\left(\mitx{B_\eg}{\hat{F}^{sr}}\right)^\frac{1}{rs(p+1)}.$$

Similarly, if $r_1>1/(p+1)$ we have the following estimate for $J_3$
$$J_3=\mitx{B_\eg}{\Fg^{\frac{1}{(p+1)}}\Fg^{\frac{p}{p+1}}|DU|^p} \le \left(\mitx{B_\eg}{\Fg^{\frac{1}{r_1^*(p+1)}}}\right)^{r_1^*}\left(\mitx{B_\eg}{F^{pr_1}}\right)^\frac{1}{r_1(p+1)}.$$

Using H\"older's inequality for the integrand $\Fg|DU|^{p+1}|Du|=F|Du|$ in $J_1$ and the above estimates, we obtain
$$\sup_\eg|\fg_\eg*J_1| \le C(\Fg,r,r_1,s)\sup \frac{|\Fg_u|}{\Fg}M(F^{s'})^\frac{1}{s'}
M(\hat{F}^{sr})^\frac{1}{rs(p+1)}
M(F^{pr_1})^\frac{1}{r_1(p+1)},$$ where \beqno{cllg}C(\Fg,r,r_1,s)=C\sup_\eg\left(\mitx{B_\eg}{\Fg^{\frac{1}{r_1^*(p+1)}}}\right)^{r_1^*}
\left(\mitx{B_\eg}{\Fg^{\frac{-s}{r^*(p+1)}}}\right)^{\frac{r^*}{s}}.\eeq 
By the definition of weights, it is clear that
\beqno{cllg1}C(\Fg,r,r_1,s)\le \left[\Fg^\frac{1}{r_1^*(p+1)}\right]_{\frac{r^*}{r_1^*s}+1}^{r_1^*}=\left[\Fg^\frac{-s}{r^*(p+1)}\right]_{\frac{r_1^*s}{r^*}+1}^{\frac{r^*}{s}}.\eeq We now choose $s,r,r_1$ such that $s'=sr=pr_1$ and $sr<2$. This is the case if $r<1$, $s=(r+1)/r$ and $r_1=(r+1)/p$. Let $$\ag(r)=\frac{1}{r_1^*(p+1)}=\frac{1}{p+1-\frac{p}{r+1}} \mbox{ and }\bg(r)=\frac{r^*}{r_1^*s}=\frac{r(p+1)-1}{r(p+1)+1}.$$

We then have $C(\Fg,r,r_1,s) \le [\Fg^{\ag(r)}]_{\bg(r)+1}^{r_1^*}$.    Clearly,  $\ag(r)$ decreases to $2/(p+2)$ and $\bg(r)$ increases to $p/(p+2)$ as $r\to1$. Thus, by the choice $\ag,\bg$ as in \mref{agbg}, if we choose $r$ close to $ 1$ then $\ag(r)<\ag$ and $\bg(r)>\bg$ and thus $[\Fg^{\ag(r)}]_{\bg(r)+1} \le C([\Fg^\ag]_{\bg+1})$, see \mref{wagbg}.  Hence, 
the above estimates give 
$$\sup_\eg|\fg_\eg*J_1| \le C([\Fg^{\ag}]_{\bg+1})\sup \frac{|\Fg_u|}{\Fg}
M(F^{sr})^\frac{2p+1}{rs(p+1)}M(\hat{F}^{sr})^\frac{1}{rs(p+1)}.$$ 

By Young's inequality we have
$$M(F^{sr})^\frac{2p+1}{rs(p+1)}M(\hat{F}^{sr})^\frac{1}{rs(p+1)}\le C(M(F^{sr})^\frac{2}{rs}+M(\hat{F}^{sr})^\frac{2}{rs}).$$

We then use \mref{HL} for $q=rs<2$ to see that
\beqno{j1est}\iidx{\Og}{\sup_\eg|\fg_\eg*J_1|} \le C([\Fg^\ag]_{\bg+1})\sup \frac{|\Fg_u|}{\Fg}
(I_1+\hat{I}_1).\eeq

Next, we write $J_2=\Fg|DU|^{p-1}|D^2U||DU|J_3$. By H\"older's inequality, $$\fg_\eg*J_2 \le \left(\mitx{B_\eg}{\Fg|DU|^{p-1}|D^2U|)^{s'}}\right)^\frac{1}{s'}
\left(\mitx{B_\eg}{\Fg^{\frac{-s}{(p+1)}}\Fg^{\frac{s}{p+1}}|DU|^s}\right)^\frac1s
J_3.$$
Again, H\"older's inequality gives
$$\left(\mitx{B_\eg}{\Fg^{\frac{-s}{(p+1)}}\Fg^{\frac{s}{p+1}}|DU|^s}\right)^\frac1s \le \left(\mitx{B_\eg}{\Fg^{\frac{-s}{r^*(p+1)}}}\right)^{\frac{r^*}{s}}\left(\mitx{B_\eg}{F^{sr}}\right)^\frac{1}{rs(p+1)}.$$

Using the estimate for $J_3$, with different $r_1$, we have
$$\sup_\eg\fg_\eg*J_2 \le C(\Fg,r,r_1,s)M((\Fg|DU|^{p-1}|D^2U|)^{s'})^\frac{1}{s'}
M(F^{sr})^\frac{1}{rs(p+1)}
M(F^{pr_1})^\frac{1}{r_1(p+1)}.$$

We will choose $s, r_1,r$ such that $pr_1=sr<2$ (and $r_1>1/(p+1)$). Again, $C(\Fg,r,r_1,s)$ can be estimated by  $C([\Fg^{\ag(r)}]_{\bg(r)+1})$, where $$\ag(r)=\frac1{p+1-\frac{p}{rs}} \mbox{ and }\bg(r)=\frac{r(p+1)-1}{rs(p+1)-p}.$$ Obviously, we choose $s>2$ and $s$ is close to 2, $r<1$ and $r$ is close to 1 such that $sr<2$, $\ag(r)<\ag$ and $\bg(r)> \bg$, see \mref{agbg}. As before, with such choice of $s,r$, we have

$$\sup_\eg\fg_\eg*J_2 \le C([\Fg^\ag]_{\bg+1})M((\Fg|DU|^{p-1}|D^2U|)^{s'})^\frac{1}{s'}M(F^{sr})^\frac{1}{rs}.$$

We have by \mref{HL} and similar estimates for \mref{dhest}, with $s_*$ is now $s'$ (which is less than 2 because $s>2$), the following
\beqno{j2est}\iidx{\Og}{\sup_\eg |\fg_\eg * J_2|} \le C([\Fg^\ag]_{\bg+1})
\sup\frac{\Fg}{\Fg_0}I_1^\frac12
I_2^\frac12.\eeq

Combining the estimates \mref{g1est}, \mref{j1est},\mref{j2est} and GN.1)  we obtain the bound \beqno{g2est}\iidx{\Og}{\sup_\eg|g_2*\fg_\eg|} \le C([\Fg^\ag]_{\bg+1})
\left[k_1I_1 + k_2I_1^\frac12 I_2^\frac12\right].\eeq

We thus prove \mref{gh1est} for the ${\cal H}^1$ norm of $g$. Therefore, by FS theorem
$$I_1 \le \|U\|_{BMO}\|g\|_{{\cal H}^1}\le C([\Fg^\ag]_{\bg+1})\|U\|_{BMO}
\left[k_1(I_1+\hat{I}_1) + k_2I_1^\frac12 I_2^\frac12\right].$$

As we noted before, see \mref{agbg}, $[\Fg^\ag]_{\bg+1}$ can be controlled by $[\Fg^\frac{2}{p+2}]_{\frac{p}{p+2}+1}$. Hence, the above gives \mref{gh1est} and the proof is complete. \eproof

\brem{GNrem1} By approximation, see \cite{SR}, the lemma also holds for $u\in W^{1,2}(\Og)$ and $U\in W^{2,2}(\Og)$ provided that the quantities $I_1,I_2$ and $\hat{I}_1$ defined in \mref{Idef} and \mref{Idef2} are finite. Furthermore, if $\Fg,\Fg_0$ are constant functions and $u=U$, the above lemma clearly gives the Gagliardo-Nirenberg inequality established in \cite{SR}.  \erem

In order to make use of the continuity in BMO norms in the asumption V) to obtain the regularity results ,  we will need the following local version of the above lemma.
 
\blemm{GNlocal} Assume as in \reflemm{GNlobal} and let $B_s,B_t$ be two concentric balls in $\Og$ with radii $t>s>0$. 
We set \beqno{Ideft} I_1(t):=\iidx{B_t}{\Fg^2(u)|DU|^{2p+2}},\,\hat{I}_1(t):=\iidx{B_t}{\Fg^2(u)|Du|^{2p+2}},\eeq and \beqno{Idef2}
I_2(t):=\iidx{B_t}{\Fg_0^2(u)|DU|^{2p-2}|D^2U|^2}.\eeq

Let $\psi$ be a $C^1$ function such that $\psi=1$ in $B_s$ and $\psi=0$ outside $B_t$. For any $\eg>0$ there is a constant $C(\eg)$ such that \beqno{GNlocalest}\barr{ll}\lefteqn{I_1(s)\le C(\eg)\|U\|_{BMO(B_t)}\sup_{x\in B_t}|D\psi(x)|^2\iidx{B_t}{|\Fg|^{2}(u)|DU|^{2p}}}\hspace{1cm}&\\& +C_\Fg
\|U\|_{BMO(B_t)}\left[(k_1+\eg)(I_1(t)+\hat{I}_1(t)) + k_2I_1^\frac12(t) I_2^\frac12(t)\right].\earr\eeq
\elemm

\bproof We revisit the proof of the previous lemma. Integrating by parts, we have
$$\iidx{\Og}{\Fg^2(u)\psi^2|DU|^{2p+2}} = -\iidx{\Og}{U\Div(\Fg^2(u)\psi^2|DU|^{2p}DU)}.$$

Again, we will show that  $g=\Div(\Fg^2\psi^2|DU|^{2p}DU)$ belongs to the Hardy space ${\cal H}^1$.
We write $g=g_1+g_2$ with $g_i=\Div V_i$, setting
$$V_1= \Fg(u) \psi|DU|^{p+1}\left(\Fg(u)\psi |DU|^{p-1}DU-\mitx{B_\eg}{\Fg(u)\psi |DU|^{p-1}DU}\right),$$ and
$$V_2= \Fg(u)\psi |DU|^{p+1}\mitx{B_\eg}{\Fg(u)\psi |DU|^{p-1}DU}.$$

In estimating $V_1$ we follow the proof of \reflemm{GNlobal} and replace $\Fg(u),\Fg_0(u)$ respectively by $\Fg(u)\psi(x)$ and $\Fg_0(u)\psi$. There will be the following extra term in estimating  $Dh$ in the right hand side of \mref{dhest} and it can be estimated as follows $$\left(\mitx{B_\eg}{\Fg^{s_*}(u)|D\psi|^{s_*}|DU|^{ps_*}}\right)^\frac1{s_*} \le \sup_{x\in B_t}|D\psi|\left(\mitx{B_\eg}{\Fg^{s_*}(u)|DU|^{ps_*}}\right)^\frac1{s_*}.$$

We then use the the following in the right hand side of \mref{g1a} (with $\Og=B_t$)
$$\barr{ll}\lefteqn{\sup|D\psi|\iidx{B_t}{\Psi_1M(\Fg^{s_*}(u)|DU|^{ps_*})^\frac{1}{s_*}} \le}\hspace{2cm}&\\& \eg\iidx{B_t}{\Psi_1^2}+ C(\eg)\sup_{x\in B_t}|D\psi|^2\iidx{B_t}{M(\Fg^{s_*}(u)|DU|^{ps_*})^\frac{2}{s_*}}.\earr$$

The last term can be bounded via \mref{HL} by $$C(\eg)\sup_{x\in B_t}|D\psi|^2\iidx{B_t}{\Fg^{2}(u)|DU|^{2p}}.$$

Using the fact that $|\psi|\le 1$ and taking $\Og$ to be $B_t$, the previous proof can go on and \mref{g1est} now becomes \beqno{g1estz}\barr{lll}\iidx{B_t}{\sup_\eg|g_1*\fg_\eg|} &\le& C
\left[\left(k_1+\eg\right)(I_1(t)+\hat{I}_1(t)) + k_2I_1^\frac12(t) I_2^\frac12(t)\right]\\&& +C(\eg)\sup|D\psi|^2\iidx{B_t}{\Fg^{2}(u)|DU|^{2p}}.\earr\eeq

Similarly, in considering $g_2=\mbox{div}V_2$, we will have an extra term $\Fg(u)|D\psi||DU|^{p+1}J_3$ in
$J_1$. We then use the following estimate  $$\sup_\eg |\fg_\eg*\Fg(u)|D\psi|\psi|DU|^{p+1}J_3 \le \sup|D\psi|M(\Fg(u)|DU|^{p+1})M(\Fg(u)|DU|^p),$$ and via Young's inequality
$$\iidx{B_t}{\sup_\eg |\fg_\eg*\Fg(u)|D\psi||DU|^{p+1}J_3} \le \eg I_1(t)+ C(\eg)\sup|D\psi|^2\iidx{B_t}{\Fg^{2}(u)|DU|^{2p}}.$$

Therefore the estimate \mref{g2est} is now \mref{g1estz} with $g_1$ being replaced by $g_2$. Combining the estimates for $g_1,g_2$ and using Young's inequality, we get  \beqno{gestz}\barr{ll}\lefteqn{\iidx{B_t}{\sup_\eg|g*\fg_\eg|}\le C(\eg)\sup|D\psi|^2\iidx{B_t}{|\Fg|^{2}|DU|^{2p}}}\hspace{1cm}&\\& +C([\Fg^\ag]_{\bg+1})
\left[\left(k_1+\eg\right)(I_1(t)+\hat{I}_1(t)) + k_2I_1^\frac12(t) I_2^\frac12(t)\right].\earr\eeq

The above gives an estimate for the ${\cal H}^1$ nowm of $g$. By FS theorem, we obtain
$$\barr{ll}\lefteqn{\iidx{B_t}{\Fg^2(u)\psi^2|DU|^{2p+2}}\le C(\eg)\|u\|_{BMO(B_t)}\sup|D\psi|^2\iidx{B_t}{\Fg^{2}(u)|DU|^{2p}}}\hspace{1cm}&\\& +C([\Fg^\ag]_{\bg+1})
\|u\|_{BMO(B_t)}\left[\left(k_1+\eg\right)(I_1(t)+\hat{I}_1(t)) + k_2I_1^\frac12(t) I_2^\frac12(t)\right].\earr$$

Since $\psi=1$ in $B_s$, the above yields \mref{GNlocalest} and the proof is complete. \eproof

Finally, the following lemma will be crucial in obtaining uniform estimates for the approximation sequence $\{u_k\}$.

\blemm{GNiter} Assume as in \reflemm{GNlocal} and let $B_\rg,B_R$ be two concentric balls in $\Og$ with radii $R>\rg>0$. Assume that there is a constant $C(\Fg,\Fg_0)$ depending on $C(\Fg,\Fg_0)$ such that the constant 
\beqno{vmo1}C(\Fg,\Fg_0):=C_\Fg
\left(\sup_u\frac{|\Fg_u|}{\Fg}+\sup_u\frac{\Fg}{\Fg_0}+1\right)\eeq is finite. Then
there is $\eg_0$ depending on $C(\Fg,\Fg_0)$ such that if $\|U\|_{BMO(B_R)}<\eg_0$ then there is a constant $C_0(\Fg,\Fg_0)$ such that \beqno{GNlocalestz}I_1(s)\le C_0(\Fg,\Fg_0)\|U\|_{BMO(B_R)}\left[\frac{1}{(t-s)^2}\iidx{B_t}{\Fg_0^{2}|DU|^{2p}}+
\hat{I}_1(t) + I_2(t)\right]\eeq for any $s,t$ such that $0<s<t<R$.
\elemm

For the proof of this lemma and later use, let us recall the following elementary iteration result (e.g., see \cite[Lemma 6.1, p.192]{Gius}). 

\blemm{Giusiter} Let $f,g,h$ be bounded nonnegative functions in the interval $[\rg,R]$ with $g,h$ being increasing. Assume that for $\rg \le s<t\le R$ we have $$f(s) \le [(t-s)^{-\ag} g(t)+h(t)]+\eg f(t)$$ with $C\ge0$, $\ag>0$ and $0\le\eg<1$. Then $$f(\rg) \le c(\ag,\eg)[(R-\rg)^{-\ag} g(R)+h(R)].$$ The constant $c(\ag,\eg)$ can be taken to be $(1-\nu)^{-\ag}(1-\nu^{-\ag}\nu_0)^{-1}$ for any $\nu$ satisfying $\nu^{-\ag}\nu_0<1$.\elemm

We are now ready to give the proof of \reflemm{GNiter}.

\bproof For any $s,t,\rg$ such that $0<\rg<s<t<R$, let $\psi$ be a cutoff function for $B_s,B_t$ with $|D\psi|\le 1/(t-s)$. Noting that $\|U\|_{BMO(B_t)} \le \|U\|_{BMO(B_R)}$. By a simple use of Young's inequality to the last product in \mref{GNlocalest} of \reflemm{GNlocal} and our assumption \mref{vmo1}, we can see easily that if $\eg,\eg_0$ are sufficiently small then for some $\nu_0\in(0,1)$ $$C(\Fg,\Fg_0)\|U\|_{BMO(B_R)} \le \nu_0,$$ and
$$I_1(s)\le \nu_0I_1(t)+ C(\Fg,\Fg_0)\|U\|_{BMO(B_R)}\left[\frac{1}{(t-s)^2}\iidx{B_t}{\Fg^{2}|DU|^{2p}} +\hat{I}_1(t)+I_2(t)\right].$$

Let $C_1(\Fg,\Fg_0)= C(\Fg,\Fg_0)\max\{1,\sup_u\frac{\Fg^2}{\Fg_0^2}\}$. The above yields
$$I_1(s)\le \nu_0I_1(t)+ C_1(\Fg,\Fg_0)\|U\|_{BMO(B_R)}\left[\frac{1}{(t-s)^2}\iidx{B_t}{\Fg_0^{2}|DU|^{2p}} +\hat{I}_1(t)+I_2(t)\right].$$

 It is clear that we can \reflemm{Giusiter} to $f=I_1$ to get $$I_1(\rg)\le C(\nu_0)C_1(\Fg,\Fg_0)\|U\|_{BMO(B_R)}\left[\frac{1}{(R-\rg)^2}\iidx{B_R}{\Fg_0^{2}|DU|^{2p}}+
 \hat{I}_1(R) + I_2(R)\right].$$
 
 The constant $C(\nu_0)$ can be taken to be $(1-\nu)^{-2}(1-\nu^{-2}\nu_0)^{-1}$ for any $\nu$ satisfying $\nu^{-2}\nu_0<1$. We can take $C(\nu_0)$  to be a fixed constant for $\nu_0\in(0,\frac12)$. Obviously, the above also holds for $\rg , R$ being replaced by $s,t$ and we proved the lemma. \eproof

\section{Proof of the main theorems}\label{proofsec}\eqnoset

\newcommand{\bint}{\dspl{\int}}

We now go back to the iterative scheme \mref{ga1} and prove our main theorems in this section.

The following lemma is the main vehicle of the proof of \reftheo{GNthm1}.

\blemm{Duplem} Assume A.1)-A.3), R)and V). Let $p\ge1$ be a number such that \beqno{pcond} \frac{2p-2}{2p}<\sup_{u\in\RR^m}\frac{\llg(u)}{\LLg(u)}.\eeq

If $R$ is sufficiently small then for any two concentric balls $B_\rg\subset B_R$ with center in $\bar{\Og}$ there is a constant $C(T)$ such that the following holds for all intergers $k\ge1$ \beqno{Dupx}\barr{ll}\lefteqn{ \sup_{t\in(0,T)}\iidx{B_\rg\cap \Og}{|Du_k|^{2p}} +\itQ{Q_\rg}{\llg(u_{k-1})|Du_k|^{2p-2}|D^2u_k|^2}\le}\hspace{.5cm}&\\& C_1(T)\itQ{Q_R}{\left[\llg(u_{1})|Du_2|^{2p-2}|D^2u_2|^2+\frac{|A_u(u_{1})|^2}{\llg(u_{1})}|Du_2|^{2p+2}\right]} \\&+ C_1(T)\frac{1}{(R-\rg)^2}\max_{1\le i\le k}\itQ{Q_R}{\llg((u_{k-1}))|Du_k|^{2p}}.\earr\eeq Here, $Q_R=(B_R\cap \Og)\times(0,T)$.
\elemm 

Before going to the proof, we recall the following elementary fact in \cite[Lemma 2.1]{sd}. 
 
 \blemm{ratiolemm} Assume the ellipticity condition A). Let $\ag$ be a number such that there is $\dg_\ag\in(0,1)$ such that $\frac{\ag}{2+\ag}=\dg_{\ag}\frac{\llg}{\LLg}$. We then have
 \beqno{D2v}AD\zeta D(\zeta |\zeta |^\ag)\ge
 \widehat{\llg}\zeta |^{\ag}|D\zeta |^2,\quad \widehat{\llg}= (1-\dg_{\ag}^2)\llg.\eeq
 
 \elemm

Furthermore, since $u_{k-1}, u_k$ are $C^2$ in $x$, we can differentiate \mref{ga1} with respect to $x$ to get
\beqno{ga2} (Du_k)_t=\Div((A(u_{k-1})D^2u_k + A_u(u_{k-1})Du_{k-1}Du_k)+Df(u_{k-1},Du_k)\quad k\ge1.\eeq

\bproof (Proof of \reflemm{Duplem}) We consider the interior case $B_\rg\subset B_R\subset \Og$ and leave the boundary case, when the center of $B_R$ is on the boundary $\partial\Og$,  to \refrem{boundaryrem} following the proof.
For any $s,t$ such that $0\le s<t \le R$ let $\psi$ be a cutoff function for $B_s,B_t$. That is, $\psi\equiv1$ in $B_s$ and $\psi\equiv0$ outside $B_t$ with $|D\psi|\le1/(t-s)$. Testing \mref{ga2} with $|Du_k|^{2p-2}Du_k\psi^2$. The assumption \mref{pcond} shows that $\ag=2p-2$ satisfies the condition of \reflemm{ratiolemm} so that we can find a positive constant $C(p)$ such that
$$\barr{ll}\lefteqn{\sup_{\tau\in(0,T)}\iidx{\Og}{|Du_k|^{2p}\psi^2}+C(p)\itQ{Q}{\llg(u_{k-1})|Du_k|^{2p-2}|D^2u_k|^2\psi^2} \le}\hspace{.2cm}&\\& -\itQ{Q}{[A_u(u_{k-1})Du_{k-1}Du_kD(|Du_k|^{2p-2}Du_k\psi^2) + Df(u_{k-1},Du_k)|Du_k|^{2p-2}Du_k\psi^2]}.\earr$$

For simplicity, we will assume in the sequel that $f\equiv0$. The presence of $f$ will be discussed in  \refrem{Frem1} after the proof.
Therefore,
$$\barr{ll}\lefteqn{\sup_{\tau\in(0,T)}\iidx{\Og}{|Du_k|^{2p}\psi^2}+C(p)\itQ{Q}{\llg(u_{k-1}|Du_k|^{2p-2}|D^2u_k|^2\psi^2} \le}\hspace{.2cm}&\\& C_0(p)\itQ{Q}{|A_u(u_{k-1})||Du_{k-1}||Du_k|^{2p-1}|D^2u_k|\psi^2+|A_u(u_{k-1})||Du_{k-1}||Du_k|^{2p}\psi|D\psi|}.\earr$$

Let $\Fg_0^2(u) = \llg(u)$ and $\Fg^2(u)= \frac{|A_u(u)|^2}{\llg(u)}$ as in A.2). Applying Young's inequality to the integrand of the first integral on the right of the above and the following $$\barr{lll}|A_u(u_{k-1})||Du_{k-1}||Du_k|^{2p}\psi|D\psi|&=&\Fg(u_{k-1})|Du_{k-1}||Du_k|^{p}\psi |\Fg_0||Du_k|^p|D\psi|\\ &\le& \Fg^2(u_{k-1})|Du_{k-1}|^2|Du_k|^{2p}\psi^2 + |D\psi|^2\Fg_0^2(u_{k-1})|Du_k|^{2p},\earr$$
we easily deduce
$$\barr{ll}\lefteqn{\sup_{\tau\in(0,T)}\iidx{\Og}{|Du_k|^{2p}\psi^2}+\itQ{Q}{\Fg_0^2(u_{k-1})|Du_k|^{2p-2}|D^2u_k|^2\psi^2} \le }\hspace{1cm}&\\&C_1\itQ{Q}{\Fg^2(u_{k-1})|Du_{k-1}|^2|Du_k|^{2p}\psi^2}+\sup|D\psi|^2\itQ{Q_t}{\Fg_0^2(u_{k-1})|Du_k|^{2p}}.\earr$$ Here, we denoted $Q_t=B_t\times(0,T)$.
 Again, a use of Young's inequality to the first integral on the right yields
$$\barr{ll}\lefteqn{\sup_{\tau\in(0,T)}\iidx{\Og}{|Du_k|^{2p}\psi^2}+\itQ{Q}{\Fg_0^2(u_{k-1})|Du_k|^{2p-2}|D^2u_k|^2\psi^2} \le }\hspace{.1cm}&\\&C_2\itQ{Q}{\Fg^2(u_{k-1})(|Du_{k-1}|^{2p+2}+|Du_k|^{2p+2})\psi^2}+C\sup|D\psi|^2\itQ{Q_t}{\Fg_0^2(u_{k-1})|Du_k|^{2p}}.\earr$$

By the choice of $\psi$, we obtain from the above the following
\beqno{ga3}\sup_{\tau\in(0,T)}\iidx{B_s}{|Du_k|^{2p}\psi^2}+H(s) \le C_2(B_0(t)+B_1(t))+C\frac{1}{(t-s)^2}G(t).\eeq
Here, for any fixed integer $k\ge1$, we set $$H(t) = \itQ{Q_t}{\Fg_0^2(u_{k-1})|Du_k|^{2p-2}|D^2u_k|^2},\,B_1(t)= \itQ{Q_t}{\Fg^2(u_{k-1})||Du_k|^{2p+2}},$$ and $$ \,B_0(t)= \itQ{Q_t}{\Fg^2(u_{k-1})||Du_{k-1}|^{2p+2}},\,G(t)= \itQ{Q_t}{\Fg_0^2(u_{k-1})||Du_k|^{2p}}.$$

We now apply \reflemm{GNiter}  for $u=u_{k-1}$ and $U=u_k$. We will see that our assumptions A.2) and A.3) imply the assumptions GN.1) and GN.2) of \reflemm{GNiter} for any $p\ge1$. Indeed, by our assumption \mref{CFg} on $\Fg_0,\Fg$ in A.2) the constants in of GN.1) are finite.
Furthermore, since $u_{k-1}$ is BMO with uniform bounded norm and the assumption A.3), $\Fg^\frac23(u_{k-1})$ belongs to the $A_\frac43$ class. As $\frac23\ge\frac{2}{p+2}$ and  $\frac43\le\frac{p}{p+2}+1$, $\Fg^\frac{2}{p+2}(u_{k-1})$ belongs to the $A_{\frac{p}{p+2}+1}$ class. Thus, 
the quantity $C(\Fg,\Fg_0)$ defined in \mref{vmo1}  is finite.
Also, our continuity assumption \mref{vmouni} on the BMO norm of $u_k$ implies the smallness of $C(\Fg,\Fg_0)\|u_k\|_{BMO(B_R)}$ if $R$ is small. Hence, for any given $\mu_1>0$ if $R=R(\mu_1)>0$ is sufficiently small then  we  have from \mref{GNlocalest} of \reflemm{GNiter} the following estimate.
\beqno{GNlocalestz}\barr{ll}\lefteqn{\iidx{B_s}{\Fg^2(u_{k-1})|Du_k|^{2p+2}}\le \mu_1\iidx{B_t}{\Fg_0^2(u_{k-1})|Du_k|^{2p-2}|D^2u_k|^2}+}\hspace{2cm}&\\&\mu_1\left[\iidx{B_t}{\Fg^2(u_{k-1})|Du_{k-1}|^{2p+2}}+\frac{1}{(t-s)^2}\iidx{B_t}{\Fg_0^2(u_{k-1})|Du_k|^{2p}}\right].\earr
\eeq

 Then \mref{GNlocalestz} and \mref{ga3} give a positive constant $C_2$, which is redefined and can depend on $k_1,k_2$, such that
\beqno{gaiter1}B_1(s) \le \mu_1[H(t)+B_{0}(t)+\frac{1}{(t-s)^2}G(t)] \quad \rg<s<t<R,\eeq
\beqno{gaiter2}H(s) \le C_2B_1(t) +C_2B_{0}(t)+\frac{1}{(t-s)^2}G(t)\quad \rg<s<t<R.\eeq

Let $t'=s+(t-s)/2$. Using \mref{gaiter1} with $s$ being $t'$ in the inequality  \mref{gaiter2} with $t$ being $t'$ and the fact that $H,B_{0},G$ are increasing, we get
\beqno{gaiter2a}H(s) \le C_2\mu_1 H(t) +C_2(\mu_1+1)B_{0}(t)+\frac{4(\mu_1+1)}{(t-s)^2}G(t)\quad \rg<s<t<R.\eeq

We can assume that $\mu_2=C_2\mu_1<1$. By \reflemm{Giusiter}, \mref{gaiter2a} yields 
$$H(\rg) \le  C_3[C_2(\mu_1+1)B_{0}(R)+\frac{4(\mu_1+1)}{(R-\rg)^2}G(R)].
$$ Here, $C_3=(1-\nu)^{-2}(1-\nu^{-2}\mu_2)^{-1}$ for any $\nu$ satisfying $\nu^{-2}\mu_2<1$. Obviously, the above also hold with $\rg, R$ replaced by $s,t$ with $\rg \le s<t \le R$. So,
\beqno{gaiter3}H(s) \le  C_3[C_2(\mu_1+1)B_{0}(t)+\frac{4(\mu_1+1)}{(t-s)^2}G(t)].\eeq

We now let $t'=(s+t)/2$ and use \mref{gaiter1} with $t$ being $t'$ and then \mref{gaiter3} with $s$ being $t'$ to see that 
$$\barr{lll}B_1(s)&\le& \mu_1[H(t') + B_{0}(t') + \frac{1}{(t'-s)^2}G(t')]\\
&\le& \mu_1[C_3[C_2(\mu_1+1)B_{0}(t)+\frac{4(\mu_1+1)}{(t-t')^2}G(t)] + B_{0}(t') + \frac{1}{(t'-s)^2}G(t')].\earr$$
Since $B_1,G$ are increasing functions, the above yields
\beqno{keyiter}B_1(s)\le \mu_3B_{0}(t)+C_4\frac{1}{(t-s)^2}G(t),\eeq where
$\mu_3=\mu_1(C_3C_2(\mu_1+1)+1)$ and $C_4=4(4\mu_1(\mu_1+1)+\mu_1)$. 

We now consider $B_0$. Applying \reflemm{GNiter} with $u=U=u_{k-1}$, so that $I_1=\hat{I}_1$, and $\Fg(u)=\Fg_0(u)$,
we see that
if $\|u_{k-1}\|_{BMO(B_R)}$, or $R$, is sufficiently small then there is a constant $C_0(\Fg,\Fg_0)$ such that for any $s,t$ satisfying $0<s<t<R$ \beqno{GNlocalestzmn}I_1(s)\le C_0(\Fg)\|u_{k-1}\|_{BMO(B_R)}\left[\frac{1}{(t-s)^2}\itQ{Q_t}{|\Fg(u_{k-1})|^{2}|Du_{k-1}|^{2p}}+ I_2(t)\right]\eeq with 
$$I_1(t) = \itQ{Q_t}{\Fg^2(u_{k-1})|Du_{k-1}|^{2p+2}}=B_0(t),$$
$$I_2(t) = \itQ{Q_t}{\Fg^2(u_{k-1})|Du_{k-1}|^{2p-2}|D^2u_{k-1}|^2}.$$ 

Going in back to the notation $\Fg_0(u)=\llg^\frac12(u)$, by \mref{UF}, we can split $Q_t$ into two disjoint sets $$Q_{(1)}:=\{(x,t)\in Q_t\,:\, \Fg(u_{k-1})(x,t) \le C(T)\Fg(u_{k-2})(x,t)\},$$
$$Q_{(2)}:=\{(x,t)\in Q_t\,:\, \Fg_0(u_{k-1})(x,t) \le C(T)\Fg_0(u_{k-2})(x,t)\}.$$

Recall that $\Fg(u)\le k_2\Fg_0(u)$ by A.2). On $Q_{(1)}$, we have $\Fg(u_{k-1})\le C(T)\Fg(u_{k-2})\le C(T)k_2\Fg_0(u_{k-2})$. Meanwhile, on $Q_{(2)}$, $\Fg(u_{k-1})\le k_2\Fg_0(u_{k-1})\le C(T)k_2\Fg_0(u_{k-2})$. Thus, $$\barr{lll}I_2(t) &=& \itQ{Q_t}{\Fg^2(u_{k-1})|Du_{k-1}|^{2p-2}|D^2u_{k-1}|^2}\\&\le& 
C(T)k_2\itQ{Q_t}{\Fg_0^2(u_{k-2})|Du_{k-1}|^{2p-2}|D^2u_{k-1}|^2}.\earr$$

Similarly, \beqno{g0}\itQ{Q_t}{|\Fg(u_{k-1})|^{2}|Du_{k-1}|^{2p}} \le C(T)k_2\itQ{Q_t}{|\Fg_0(u_{k-2})|^{2}|Du_{k-1}|^{2p}}.\eeq

Using these estimates in \mref{GNlocalestzmn}, we obtain $$B_0(s) \le C_1(\Fg,\Fg_0,T)\|u_{k-1}\|_{BMO(B_R)}[\frac{1}{(t-s)^2}G_0(t)+H_0(t)],$$ where $$G_0(t)=\itQ{Q_t}{|\Fg_0(u_{k-2})|^{2}|Du_{k-1}|^{2p}},$$ $$H_0(t) =\itQ{Q_t}{\Fg_0^2(u_{k-2})|Du_{k-1}|^{2p-2}|D^2u_{k-1}|^2}.$$

Using the above estimate for $B_0$ in \mref{gaiter3} and \mref{keyiter} and adding the results, we can easily see that if $\|u_{k-1}\|_{BMO(B_R)}$ is sufficiently small then \beqno{keyitermmm}H(s)+B_1(s) \le \mu_4 H_0(t) + \frac{C_5}{(t-s)^2}[G(t)+G_0(t)],\eeq for some  $C_5$ depends on $\Fg,\Fg_0,k_1,k_2,T$ and $$\mu_4= C_1(\Fg,\Fg_0,T)\|u_{k-1}\|_{BMO(B_R)}[C_3[C_2(\mu_1+1)+\mu_3].$$

We now define
$${\cal B}_k(t)=\itQ{Q_t}{[\llg(u_{k-1})|Du_k|^{2p-2}|D^2u_k|^2+\Fg^2(u_{k-1})||Du_{k}|^{2p+2}]},$$
$${\cal G}_k(t)=\itQ{Q_t}{[\llg(u_{k-1})|Du_k|^{2p}+\llg(u_{k-2})|Du_{k-1}|^{2p}]}.$$

We then have from \mref{keyitermmm} that
\beqno{Biter2}{\cal B}_k(s) \le \mu_4{\cal B}_{k-1}(t) + \frac{C_5}{(t-s)^2}{\cal G}_k(t).\eeq

As before, we can assume that $R$ is sufficiently small such that $\mu_4<1$. For any $a\in(0,1)$ such that $\mu_4a^{-2}<1$ we define the sequences $t_0=\rg$ and $t_{i+1}=t_i+(1-a)a^i(R-\rg)$. Iterate the above $k-2$ times to get
$$\barr{lll}{\cal B}_k(\rg) &\le& \mu_4^k{\cal B}_2(t_{k-2})+\sum_{i=0}^{k-2}\mu_4^ia^{-2i}\frac{C_5}{(1-a)^2(R-\rg)^2}{\cal G}_{k-i}(t_{i+1})\\&\le& {\cal B}_2(R)+\frac{C_6(a,\mu_4)}{(R-\rg)^2}\max_{2\le i\le k}{\cal G}_i(R).\earr$$

This shows that the quantity $$\itQ{Q_t}{[\llg(u_{k-1})|Du_k|^{2p-2}|D^2u_k|^2+\Fg^2(u_{k-1})||Du_{k}|^{2p+2}]}, \,t\ge\rg,\, k\ge1$$ can be bounded by
$$ C_1(T)\itQ{Q_R}{\left[\llg(u_{1})|Du_2|^{2p-2}|D^2u_1|^2+\frac{|A_u(u_{1})|^2}{\llg(u_{1})}|Du_2|^{2p+2} + \frac{1}{(R-\rg)^2}\max_{1\le i\le k}\llg((u_{k-1}))|Du_k|^{2p}\right]}.$$ 

Using this and \mref{g0} in \mref{ga3}, and the estimate for $B_0(t)$, we obtain \mref{Dupx} of the lemma. \eproof

\brem{Frem1} If $f\ne0$ then the growth assumption F) gives
\beqno{fcond2} |f(u_{k-1},Du_k)| \le C|Du_k| + C|u_{k-1}|^b+C.\eeq
\beqno{fcond3} |Df(u_{k-1},Du_k)| \le C|D^2 u_k| + C|u_{k-1}|^{b-1}|Du_{k-1}|+C.\eeq

Testing \mref{ga2} with $|Du_k|^{2p-2}Du_k\psi^2$, we will have the extra term $Df(u_{k-1},Du_k)|Du_k|^{2p-1}\psi^2$ on the right of our estimates in the proof. For any positive $\eg>0$ we can use Young's inequality to have  $$\barr{lll}Df(u_{k-1},Du_k)|Du_k|^{2p-1}&\le&C(|D^2u_k|+ |u_{k-1}|^{b-1}|Du_{k-1}|+C)|Du_k|^{2p-1}\\&\le& \eg|Du_k|^{2p-2}|D^2u_k|^2+ \\ && C(p,\eg)[\Fg^{-\frac{k_p}{p+1}}|u_{k-1}|^{k_p(b-1)}+\Fg^2|Du_{k-1}|^{2p+2}+|Du_k|^{2p}],\earr$$ where $\frac{1}{k_p} =1- \frac{2p-1}{2p}-\frac{1}{2p+2}$. Since $u_k, u_{k-1}$ are BMO and $\Fg$ is bounded from below, the integral of the first term on the right of the above inequalities is bounded.
\mref{gaiter2} now becomes
\beqno{gaiter2z}H(s) \le \eg H(t)+ C_2B_1(t) +C_2B_0(t)+\frac{1}{(t-s)^2}G(t)+C(\eg)\quad \rg<s<t<R.\eeq

Choosing $\eg$ small, we can see that the iteration arguments in the proof are still in force and the proof can continue.
\erem

\brem{boundaryrem} We discuss the case when the centers of $B_\rg, B_R$ are on the boundary $\partial \Og$. Let us first consider the case  $u$ satisfies the Dirichlet condition $u=0$ on $\partial \Og$. By flattening the boundary we can assume that $B_R\cap\Og$ is the set $$B^+=\{x\,:\, x=(x_1,\ldots,x_n)\mbox{ with } x_n\ge0 \mbox{ and } |x|<R\}.$$

For any point $x=(x_1,\ldots,x_n)$ we denote by $\bar{x}$ its reflection  across the plane $x_n=0$, i.e.,  $\bar{x}=(x_1,\ldots,-x_n)$. Accordingly, we denote by $B^-$ the reflection of $B^+$. For $u=u_k$ we define the odd reflection of $u$ by $\bar{u}$, i.e. $\bar{u}(x,t)=-u(\bar{x},t)$ for $x\in B^-$. We then consider the odd extension $U$ in $B=B^+\cup B^-$ of $u$  
$$U(x,t)=\left\{\barr{ll}u(x,t) &\mbox{if $x\in B^+$},\\\bar{u}(x,t)&\mbox{if $x\in B^-$}.\earr \right.$$
It is easy to see that $\bar{u}$ satisfies in $B^-$ a system similar to \mref{ga2} for $u$ in $B^+$. As in the proof of the lemma, we test the system for $u_k$ with $|Du_k|^{2p-2}Du_k\psi^2$ and the system for $\bar{u}$ with $|D\bar{u}_k|^{2p-2}D\bar{u}_k\psi^2$ and then sum the results. The integration parts results the extra boundary terms along the flat boundary parts $\partial B^+$ and $\partial B^-$. Using the facts that either $D_{x_i}u=D_{x_i}\bar{u}=0$ for $i\ne n$ or $D_{x_n}u=D_{x_n}\bar{u}$ and the outward normal vectors of $B^+$ and $B^-$ are opposite we can easily see that those boundary terms are either zero or cancel each others in the summation. Thus, we can obtain \mref{ga3} again with $u_{k-1},u_k$ being replaced by $U_{k-1},U_k$. Since $U_k$ belong to $W^{2,\infty}(B)$ the argument can continue and the lemma holds for $U_k$ and then $u_k$.

The same argument applies for the Neumann boundary condition if we we use the even extension for $u_k$.

\erem

We now give the proof of \reftheo{GNthm1}.

\bproof  We test the systems \mref{ga1} with $u_k$ and use Young's inequality to have
$$\sup_{t\in[0,T]}\iidx{\Og}{|u|^{2}}+\itQ{Q}{\llg(u_{k-1})|Du_k|^{2}}\le C\itQ{Q}{[|u_{k-1}|^2+|u_{k-1}|^{b+1}+1]}.$$

The uniform bound assumption on the BMO norms of $u_{k-1}$ yields that the right hand side is bounded uniformly for all $k$. Thus, there is a constant $C$ such that $$\itQ{Q}{\llg(u_{k-1})|Du_k|^2} \le C\quad \forall k.$$

Now, for any $0<\rg<R$ and concentric balls $B_\rg ,B_R$ with centers in $\bar{\Og}$ let us assume that there is some $p\ge1$ such that there is a constant $C_0(\rg,R,u_0)$ depending on $\rg,R$ and $\sup_{t\in(0,T)}\|u_0(\cdot,t)\|_{C^1(\Og)}$ on such that ($Q_R=B_R\times(0,T)$) \beqno{startiter}\itQ{Q_\rg}{\llg(u_{k-1})|Du_k|^{2p}} \le C_0(\rg,R,u_0)\quad \forall k.\eeq

It is well known that the $C^1$ norms of $u_1$ and $u_2$ can be bounded by that of $u_0$. Now, if $p$ satisfies \mref{pcond} then \reflemm{Duplem} and \mref{startiter} establish the existence of a constant $C_1(\rg,R)$ such that the following holds for all integers $k$
\beqno{Viter1} \sup_{t\in(0,T)}\iidx{B_\rg}{|Du_k|^{2p}} +\itQ{Q_\rg }{\Fg_0^2(u_{k-1})|Du_k|^{2p-2}|D^2u_k|^2}\le C_1(\rg,R,u_0)\eeq if $0<\rg<R$ and $R$ is sufficiently small.

Let $\chi_0$ be any number such that $1<\chi_0<1+\frac2n$. Denote $V=|Du_k|^{p}$ and use H\"older's inequality to get $$\itQ{Q}{\llg V^{2\chi_0}} \le \left(\itQ{Q}{\llg^{r}}\right)^\frac1r\left(\itQ{Q}{V^{2(1+\frac{2}{n})}}\right)^\frac1{r'}$$ where $\llg=\llg(u_{k-1})$ and $r$ is a number such that $r'\chi_0=1+\frac2n$. 

Recall the Sobolev imbedding inequality $$\|V\|_{L^{\frac{2(n+2)}{n}}(Q)} \le C\sup_t\|V(\cdot,t)\|_{L^2(\Og)} +C\left(\itQ{Q}{|DV|^2}\right)^\frac12$$ and the fact that $u_k$ is BMO so that $\llg(u_{k-1})$ belongs to $L^r(\Og)$ for any $r>1$ (see \cite{Gius}). The above estimates for $Q=Q_\rg$ show that there is a constant $C(\rg)$ such that \beqno{Viter}\itQ{Q_\rg}{\llg V^{2\chi_0}} \le C(\rg)\left[\sup_t\|V(\cdot,t)\|_{L^2(B_\rg)} +\left(\itQ{Q_\rg}{|DV|^2}\right)^\frac12\right].\eeq

From the ellipticity condition A) and \mref{Viter1} we see that  the right hand side is bounded. Hence 
$$\itQ{Q}{\llg(u_{k-1})|Du_k|^{2p\chi_0}} \le C_2(\rg,R,u_0)\quad \forall k.$$

Therefore, \mref{startiter} holds again with $p$ now is $p\chi_0$. We already showed that \mref{startiter} is valid for $p=1$. Thus, we can repeat the argument $k$ times until $2\chi_0^k>n$ as long as the ratio condition \mref{pcond} of \reflemm{Duplem} is verified for $p=\chi_0^k$.
The assumption R) shows that we can choose $\chi_0,k$ such that $1<\chi_0<1+\frac2n$, $2\chi_0^k>n$ and the ratio condition \mref{pcond} holds for $p=\chi_0^k$. Therefore,  \mref{Viter1} holds for $2p=2\chi_0^k$ . We now cover $\Og$ with finitely many balls of radius $R/2$ to obtain
\beqno{Viter1z} \sup_{t\in(0,T)}\iidx{\Og}{|Du_k|^{2p}} +\itQ{Q}{\Fg_0^2(u_{k-1})|Du_k|^{2p-2}|D^2u_k|^2}\le C(\Og,T,u_0).\eeq

For each $t\in(0,T)$, \mref{Viter1z} shows that the norms $\|u_k(\cdot,t)\|_{W^{1,2p}(\Og)}$ for some $2p>n$ are bounded uniformly in $t$ by a constant depending only on the size of $\Og,T$ and $u_0$. By Sobolev's imbedding theorem $\{u_k(\cdot,t)\}$ is a bounded sequence in $C^{\ag}(\Og)$ for some $\ag>0$. From the system for $u_k$, \mref{Viter1z} with $p=1$ also shows that $\|(u_k)_t\|_{L^2(Q)}$ is uniformly bounded. Together, these facts show that the solutions $u_k$ are uniformly H\"older continuous in $(x,t)$ and that $\{u_k\}$ is  bounded in $C^{\bg}(Q)$ for some $\bg>0$ We then see that there is a relabeled subsequence $\{u_k\}$ converges in $C^0(Q)$ to some $u$. Using difference quotient in $t$ we see that $u_t\in{L^2(Q)}$. The above estimate \mref{Viter1z} also shows that we can assume $Du_{k+1}(\cdot,t)$ converges weakly to $Du(\cdot,t)$ in $L^2(\Og)$ for each $t\in(0,T))$. By the continuity of $A$ in its variable $u$, we see that $u$ weakly solves \mref{e1}.

By the semicontinuity of norms, \mref{Viter1z} implies  \beqno{Viter1zz} \sup_{t\in(0,T)}\iidx{\Og}{|Du|^{2p}}+\itQ{Q}{\llg(u)|Du|^{2p-2}|D^2u|^2} \le C(\Og,u_0).\eeq

Since $2p>n$, the above implies that $u$ is H\"older continuous and its regularity in $x$. Since $u_t$ is in $L^2(Q)$. It is easy to derive from these facts that $u$ is H\"older in $(x,t)$. By \cite{GiaS}, $Du$ is H\"older in $(x,t)$ and then $u$ is a classical solution.  \eproof

We now turn to the proof of our second theorem. 

\bproof {\em (Proof of \reftheo{Amthm1})}
Let $u$ be the classical solution of the system \mref{e1} in $\Og\times(0,T_0)$. We can differentiate \mref{e1} to have
\beqno{ga4} (Du)_t=\Div((A(u)D^2u + A_u(u)DuDu)+Df(u,Du).\eeq

For any $s,t$ such that $0<s<t<R$ let $\psi$ be a cutoff function for two concentric balls $B_s,B_t$ with centers in $\bar{\Og}$. That is, $\psi\equiv1$ in $B_s$ and $\psi\equiv0$ outside $B_t$ with $|D\psi|\le1/(t-s)$. As in the proof of the previous theorem, we test \mref{ga4} with $|Du|^{2p-2}Du\psi^2$. Since $u\in C^{1,2}(\bar{\Og}\times(0,T_0)$, the local Gagliardo-Nirenberg inequality of \reflemm{GNlocal} applies here for $U=u$. Thus, as in the proof of \reflemm{Duplem}, we can use \reflemm{GNiter} to obtain \beqno{GNlocalestz1}\barr{ll}\lefteqn{\iidx{B_s}{\Fg^2|Du|^{2p+2}}\le \mu_1\iidx{B_t}{\Fg_0^2|Du|^{2p-2}|D^2u|^2}+}\hspace{2cm}&\\&\mu_1\left[\iidx{B_t}{\Fg^2|Du|^{2p+2}}+\frac{1}{(t-s)^2}\iidx{B_t}{\Fg^2|Du|^{2p}}\right].\earr
\eeq Here, $\Fg_0^2 = \llg(u)$ and $\Fg^2= \frac{A_u^2(u)}{\llg(u)}$. We now set $$H(t) = \itQ{Q_t}{\Fg_0^2(u)|Du|^{2p-2}|D^2u|^2},\,B(t)= \itQ{Q_t}{\Fg^2(u)||Du|^{2p+2}},$$ $$ G(t)= \itQ{Q_t}{\Fg^2(u)||Du|^{2p}}.$$

Because A.4) obviously holds for $u_k=u_{k-1}=u$, it is clear that the proof of \reflemm{Duplem} with $H,H,B_0,B_1,{\cal B}_k$, and $G$ being replaced by the new definitions,  and $B_0=B_{1}={\cal B}_k=B$, now leads to

\beqno{keyiter1}B(s)\le \mu_4B(t)+C_4\frac{1}{(t-s)^2}G(t),\quad 0<\rg<s<t<R,\eeq where
$\mu_4<1$.  For any $a\in(0,1)$ such that $\mu_4a^{-2}<1$ we define the sequences $t_0=\rg$ and $t_{i+1}=t_i+(1-a)a^i(R-\rg)$. Iterate the above to get
$$B(\rg) \le \mu_4^kB(t_k)+\sum_{i=0}^{k-1}\mu_4^ia^{-2i}\frac{C_4}{(1-a)^2(R-\rg)^2}G(t_k).$$ Let $k$ tend to infinity and use the fact that $\mu_4\in(0,1)$ and $B(R)$ is finite to get
$$B(\rg) \le\frac{C_5(a,\mu_4)}{(R-\rg)^2}G(R).$$

We now see that a similar argument in the proof of \reftheo{GNthm1} with $u_k$ being $u$ now gives \beqno{Viter11} \sup_{t\in(0,T)}\iidx{B_\rg}{|Du|^{2p}} +\itQ{Q_\rg }{\Fg_0^2|Du|^{2p-2}|D^2u|^2}\le C_1(\rg,R)\eeq if $0<\rg<R$ and $R$ is sufficiently small and some $p$ such that $2p>n$. Finite covering $\Og$ with balls $B_{R/2} $ yields
\beqno{Viter1zz1} \sup_{t\in(0,T)}\iidx{\Og}{|Du|^{2p}}+\itQ{Q}{\llg(u)|Du|^{2p-2}|D^2u|^2} \le C(\Og,R).\eeq

Hence $u$ is H\"older continuous and its regularity in $x$. From the system for $u$ and the above, with $p=1$, we see that $u_t$ is in $L^2(Q)$. It is now standard to show that $u$ is H\"older in $(x,t)$ and $Du$ is H\"older continuous. We now can refer to Amann's results to see that $u$ exists globally.  \eproof

\bibliographystyle{plain}

\begin{thebibliography}{10}

\bibitem{Am1} H. Amann. \newblock{Dynamic theory of quasilinear parabolic equations II. Reaction–diffusion systems,} {\em Differential Integral
Equations,} Vol. 3, no. 1 (1990), pp. 13–-75. 

\bibitem{Am2} H. Amann. \newblock{Dynamic theory of quasilinear parabolic systems III. Global existence,} {\em  Math Z.} 202 (1989), pp. 219–-250.

\bibitem{sd} S. Ahmad and D. Le. \newblock Global and Blow Up Solutions to Cross Diffusion Systems.
\newblock{Nonlinear Analysis Series A: TMA}. In press.


\bibitem{AF} A. Friedman. \newblock{\em Partial Differential Equations,} New York, 1969.

\bibitem{FPW} B. Franchi, C. Perez and R. L. Wheeden. \newblock Self-Improving Properties of John Nirenberg and
Poincar\'e Inequalities on Spaces of Homogeneous Type. \newblock{\em J.  Functional Analysis}, 153, 108--146, 1998.

\bibitem{Gius}
E. Giusti.
\newblock {\em Direct Methods in the Calculus of Variations}.
\newblock World Scientific, 2003.

\bibitem{JN} R. L. Johnson and C. J. Neugebauer. \newblock Properties of BMO functions whose reciprocals are also BMO. \newblock{\em Z. Anal. Anwendungen}, 12(1):3-11, 1993.


\bibitem{letrans} D. Le. \newblock Regularity of BMO weak solutions to nonlinear parabolic systems via homotopy. \newblock{Trans. Amer. Math. Soc}. 365 (2013), no. 5, 2723--2753.

\bibitem{leans} D. Le. \newblock Global existence results for near triangular nonlinear parabolic systems. \newblock {Adv. Nonlinear Studies}. 13 (2013), no. 4, 933-944.

\bibitem{DLT}  D. Le, L. Nguyen and T. Nguyen. \newblock{Coexistence in Cross Diffusion systems.} {\em Indiana Univ. J. Math.} Vol.56, No. 4, pp.1749-1791, 2007.



\bibitem{OP}
J. Orobitg and C. P\'erez.
\newblock $A_p$ weights for nondoubling measures in $\RR^n$ and applications.
newblock {\em 
Transactions of the American mathematical society}, 354 (2002), 2013-2033.

\bibitem{SR} P. Strzelecki.  \newblock{Gagliardo Nirenberg inequalities with a BMO term.} {\em Bull. London Math. Soc.} Vol. 38, pp. 294-300, 2006.

\bibitem{SKT} N. Shigesada, K. Kawasaki  and E.~Teramoto. \emph{ 
Spatial segregation of interacting species}.  J. Theor. Biol., 
79(1979), 83-- 99.

\bibitem{st} E. M. Stein. \newblock {\em Harmonic Analysis, Real Variable Methods, Orthogonality and Oscillatory Integrals.}
Princeton Univ. Press, Princeton, NJ, 1993.

\end{thebibliography}

\end{document}